\documentclass[11pt]{amsart}

\author[C.~Sanna]{Carlo Sanna$^\dagger$}
\thanks{$^\dagger$C.~Sanna is a member of the INdAM group GNSAGA}
\address{Universit\`a degli Studi di Torino\\Department of Mathematics\\Turin, Italy}
\email{carlo.sanna.dev@gmail.com}

\title[On numbers divisible by the product of their nonzero base $b$ digits]{On numbers divisible by the product of \\their nonzero base $b$ digits}

\usepackage{amsmath}
\usepackage{amssymb}
\usepackage{amsthm}
\usepackage{geometry}
\usepackage{color}
\usepackage{hyperref}
\usepackage{enumerate}
\hypersetup{colorlinks=true}
\usepackage{placeins}

\newtheorem{thm}{Theorem}[section]

\begin{document}

\begin{abstract}
For each integer $b \geq 3$ and every $x \geq 1$, let $\mathcal{N}_{b,0}(x)$ be the set of positive integers $n \leq x$ which are divisible by the product of their nonzero base $b$ digits.
We prove bounds of the form $x^{\rho_{b,0} + o(1)} < \#\mathcal{N}_{b,0}(x) < x^{\eta_{b,0} + o(1)}$, as $x \to +\infty$, where $\rho_{b,0}$ and $\eta_{b,0}$ are constants in ${]0,1[}$ depending only on $b$.
In particular, we show that $x^{0.526} < \#\mathcal{N}_{10,0}(x) < x^{0.787}$, for all sufficiently large $x$.
This improves the bounds $x^{0.495} < \#\mathcal{N}_{10,0}(x) < x^{0.901}$, which were proved by De~Koninck and Luca.
\end{abstract}

\maketitle

\section{Introduction}

Let $b \geq 2$ be an integer. 
Then, every positive integer $n$ has a unique representation as
\begin{equation*}
n = \sum_{j\,=\,0}^\ell d_j b^j, \quad d_0, \dots, d_\ell \in \{0, \dots, b - 1\}, \quad d_\ell \neq 0,
\end{equation*}
where $d_0, \dots, d_\ell$ are the \emph{base $b$ digits} of $n$.
Positive integers whose base $b$ digits obey certain restrictions have been investigated by several authors.
For instance, an asymptotic formula for the counting function of \emph{$b$-Niven numbers}, that is, positive integers divisible by the sum of their base $b$ digits, has been proved by De~Koninck, Doyon, and K\'atai~\cite{MR1957109}, and (independently) by Mauduit, Pomerance, and S\'ark\"ozy~\cite{MR2166377}.
Also, arithmetic properties of integers with a fixed sum of their base $b$ digits have been studied by Luca~\cite{MR2294785}, Mauduit and S\'ark\"ozy~\cite{MR1456239}.
Moreover, prime numbers with specific restrictions on their base $b$ digits have been investigated by Bourgain~\cite{MR3047097, MR3319636} and Maynard~\cite{Maynard2, Maynard1} (see~\cite{MR1752255, MR3710628} for similar works on almost primes and squarefree numbers).

Let $p_b(n)$ be the product of the base $b$ digits of $n$, and let $p_{b,0}(n)$ be the product of the nonzero base $b$ digits of $n$.
For all $x \geq 1$, define the sets
\begin{equation*}
\mathcal{N}_b(x) := \big\{n \leq x : p_b(n) \mid n\big\} \quad\text{ and }\quad \mathcal{N}_{b,0}(x) := \big\{n \leq x : p_{b,0}(n) \mid n\big\} .
\end{equation*}
Note that $\mathcal{N}_b(x) \subseteq \mathcal{N}_{b,0}(x)$ and that $n \in \mathcal{N}_b(x)$ implies that all the base $b$ digits of $n$ are nonzero.
Furthermore, $\mathcal{N}_2(x) = \big\{2^k - 1 : k \geq 1\big\}$ and $\mathcal{N}_{2,0}(x) = \mathbb{N}$.
Hence, in what follows, we will focus only on the case $b \geq 3$.

De~Koninck and Luca~\cite{MR2298113} (see also~\cite{MR3734412} for the correction of a numerical error in~\cite{MR2298113}) studied $\mathcal{N}_{10}(x)$ and $\mathcal{N}_{10,0}(x)$.
They proved the following bounds.

\begin{thm}\label{thm:konluc}
We have
\begin{equation*}
x^{0.122} < \#\mathcal{N}_{10}(x) < x^{0.863}
\end{equation*}
and
\begin{equation*}
x^{0.495} < \#\mathcal{N}_{10,0}(x) < x^{0.901}
\end{equation*}
for all sufficiently large $x$.
\end{thm}

In this paper, we prove some bounds for the cardinalities of $\mathcal{N}_b(x)$ and $\mathcal{N}_{b,0}(x)$.
In particular, for $b = 10$, we get the following improvement of three of the bounds of Theorem~\ref{thm:konluc}.

\begin{thm}\label{thm:base10}
We have
\begin{equation*}
\#\mathcal{N}_{10}(x) < x^{0.717}
\end{equation*}
and
\begin{equation*}
x^{0.526} < \#\mathcal{N}_{10,0}(x) < x^{0.787}
\end{equation*}
for all sufficiently large $x$.
\end{thm}

\subsection*{Notation}
We use the Landau--Bachmann ``little oh'' notation $o$, as well as the Vinogradov symbol $\ll$.
We omit the dependence on $b$ of the implied constants.
We write $P(n)$ for the greatest prime factor of an integer $n > 1$.
As usual, $\pi(x)$ denotes the number of prime numbers not exceeding $x$.
We write $\nu_p$ for the $p$-adic valuation.

\section{Upper bounds}

For every $s \geq 0$, let us define
\begin{equation*}
\zeta_b(s) := \sum_{d \,=\, 1}^{b - 1} \frac1{d^s} .
\end{equation*}
We give the following upper bounds for $\#\mathcal{N}_{b,0}(x)$ and $\#\mathcal{N}_b(x)$.

\begin{thm}\label{thm:Nb0upper}
Let $b \geq 3$ be an integer.
We have
\begin{equation*}
\#\mathcal{N}_{b,0}(x) < x^{\eta_{b,0} + o(1)} ,
\end{equation*}
as $x \to +\infty$, where
\begin{equation*}
\eta_{b,0} := 1 + \frac1{(1 + s_{b,0})\log b}\log\!\left(\frac{1 + \zeta_b(s_{b,0})}{b}\right) \in {]0, 1[}
\end{equation*}
and $s_{b,0}$ is the unique solution of the equation
\begin{equation}\label{equ:Equb0}
\frac{(1 + s)\zeta_b^\prime(s)}{1 + \zeta_b(s)} - \log\!\left(\frac{1 + \zeta_b(s)}{b}\right) = 0
\end{equation}
over the positive real numbers.
\end{thm}

\begin{thm}\label{thm:Nbupper}
Let $b \geq 3$ be an integer.
We have
\begin{equation*}
\#\mathcal{N}_b(x) < x^{\eta_b + o(1)} ,
\end{equation*}
as $x \to +\infty$, where $\eta_3 := \log 2 / \log 3$,
\begin{equation*}
\eta_b := 1 + \frac1{(1 + s_b)\log b}\log\!\left(\frac{\zeta_b(s_b)}{b}\right), \quad b \geq 4,
\end{equation*}
and $s_b$ is the unique solution of the equation
\begin{equation}\label{equ:Equb}
\frac{(1 + s) \zeta_b^\prime(s)}{\zeta_b(s)} - \log\!\left(\frac{\zeta_b(s)}{b}\right) = 0
\end{equation}
over the positive real numbers.
\end{thm}

We remark that for $b = 3$ the bound of Theorem~\ref{thm:Nbupper} is obvious.
Indeed, it is an easy consequence of the fact that all the base $3$ digits of each $n \in \mathcal{N}_3(x)$ are equal to $1$ or $2$.
We included it just for completeness.

Using the PARI/GP~\cite{PARI2} computer algebra system, the author computed $s_{10,0} = 1.286985\!\dots$ and $s_{10} = 1.392189\!\dots$, which in turn give $\eta_{10,0} = 0.7869364\!\dots$ and $\eta_{10} = 0.7167170\!\dots$
Hence, the upper bounds of Theorem~\ref{thm:base10} follow.

\subsection*{Proof of Theorem~\ref{thm:Nb0upper}}

First, we shall prove that Equation~\eqref{equ:Equb0} has a unique positive solution.
For $s \geq 0$, let
\begin{equation*}
F_b(s) := \frac{(1 + s)\zeta_b^\prime(s)}{1 + \zeta_b(s)} - \log\!\left(\frac{1 + \zeta_b(s)}{b}\right) .
\end{equation*}
Since $b \geq 3$, we have
\begin{equation}\label{equ:Fb0oo}
F_b(0) = -\frac{\log((b-1)!)}{b} < 0 \quad\text{ and }\quad \lim_{s \,\to\, +\infty} F_b(s) = \log\!\left(\frac{b}{2}\right) > 0 .
\end{equation}
Furthermore, a bit of computation shows that
\begin{equation}\label{equ:Fbprime}
F_b^\prime(s) = \frac{(1 + s)\left((1 + \zeta_b(s))\zeta_b^{\prime\prime}(s) - (\zeta_b^\prime(s))^2\right)}{(1 + \zeta_b(s))^2} > 0 ,
\end{equation}
for all $s \geq 0$, since, by Cauchy--Schwarz inequality, we have
\begin{equation}\label{equ:CS}
(\zeta_b^\prime(s))^2 = \left(-\sum_{d \,=\, 1}^{b-1} (\log d) d^{-s}\right)^2 < \left(\sum_{d \,=\, 1}^{b-1} d^{-s}\right)\left(\sum_{d \,=\, 1}^{b-1} (\log d)^2 d^{-s}\right) = \zeta_b(s)\zeta_b^{\prime\prime}(s) .
\end{equation}
At this point, by \eqref{equ:Fb0oo} and \eqref{equ:Fbprime}, it follows that Equation~\eqref{equ:Equb0} has a unique positive solution.

Let us assume $x \geq 1$ sufficiently large, and let $\alpha \in {]0, 1[}$ be a constant (depending on $b$) to be determined later.
Also, let $P_b$ be the greatest prime number less than $b$, and define the set
\begin{equation*}
\mathcal{N}_b^{\,\prime}(x) := \big\{n \leq x : d \mid n \text{ for some } d > x^\alpha \text{ with } P(d) \leq P_b \big\} .
\end{equation*}
Suppose $n \in \mathcal{N}_b^{\,\prime}(x)$.
Then there exists $d > x^\alpha$ with $P(d) \leq P_b$ such that $d \mid n$.
Clearly, for any fixed $d$, there are at most $x / d$ possible values for $n$.
Moreover, setting
\begin{equation*}
\mathcal{S}(t) := \big\{d \leq t : P(d) \leq P_b \big\} ,
\end{equation*}
it follows easily that $\#\mathcal{S}(t) \ll (\log t)^{\pi(P_b)}$ for all $t > 2$.
Therefore, we have
\begin{align*}
\#\mathcal{N}_b^{\,\prime}(x) &\leq \sum_{\substack{x^\alpha \,<\, d \,\leq\, x}} \frac{x}{d} = x\left(\left.\frac{\#\mathcal{S}(t)}{t}\right|_{t\,=\,x^\alpha}^x + \int_{t\,=\,x^\alpha}^x \frac{\#\mathcal{S}(t)}{t^2}\,\mathrm{d}t\right) \ll (\log x)^{\pi(P_b)}\left(1 + x^{1 - \alpha}\right) ,
\end{align*}
and consequently
\begin{equation}\label{equ:N1upperbound}
\#\mathcal{N}_b^{\,\prime}(x) < x^{1 - \alpha + o(1)} ,
\end{equation}
as $x \to +\infty$.

Now suppose $n \in \mathcal{N}_{b,0}^{\,\prime\prime}(x) := \mathcal{N}_{b, 0}(x) \setminus \mathcal{N}_b^{\,\prime}(x)$.
Put $N := \lfloor \log x / \log b \rfloor + 1$, so that $n$ has at most $N$ base $b$ digits.
For each $d \in \{1, \dots, b - 1\}$, let $N_d$ be the number of base $b$ digits of $n$ which are equal to $d$.
Also, let $N_0 := N - (N_1 + \cdots + N_{b-1})$.
Hence, $N_0, \dots, N_{b-1}$ are nonnegative integers such that $N_0 + \cdots + N_{b - 1} = N$.
Furthermore,
\begin{equation*}
p_{b,0}(n) = 1^{N_1} \cdots (b-1)^{N_{b-1}} \leq x^{\alpha} < b^{\alpha N} .
\end{equation*}
Let $\beta > 0$ be a constant (depending on $b$) to be determined later.
For fixed $N_0, \dots, N_{b-1}$, by elementary combinatorics, the number of possible values for $n$ is at most
\begin{equation*}
\frac{N!}{N_0! \cdots N_{b-1}!} .
\end{equation*}
Hence, summing over all possible values for $N_0, \dots, N_{b - 1}$, we get
\begin{align*}
\#\mathcal{N}_{b,0}^{\,\prime\prime}(x) &\leq \sum_{\substack{N_0 + \cdots + N_{b - 1} \,=\, N \\ 1^{N_1} \cdots (b-1)^{N_{b - 1}} \,\leq\, b^{\alpha N}}} \frac{N!}{N_0! \cdots N_{b-1}!} \\
&\leq \sum_{N_0 + \cdots + N_{b - 1} \,=\, N} \frac{N!}{N_0! \cdots N_{b-1}!} \left(\frac{b^{\alpha N}}{1^{N_1} \cdots (b-1)^{N_{b-1}}}\right)^\beta \\
&= \left(b^{\alpha\beta} (1 + \zeta_b(\beta))\right)^N ,
\end{align*}
where we employed the multinomial theorem.
Therefore, since $N \leq \log x / \log b + 1$, we have
\begin{equation}\label{equ:N2upperbound}
\#\mathcal{N}_{b,0}^{\,\prime\prime}(x) < x^{\gamma + o(1)} ,
\end{equation}
as $x \to +\infty$, where 
\begin{equation}\label{equ:gammaalphabeta}
\gamma := \alpha\beta + \frac{\log(1 + \zeta_b(\beta))}{\log b} .
\end{equation}

At this point, in light of \eqref{equ:N1upperbound} and \eqref{equ:N2upperbound}, we shall choose $\alpha$ and $\beta$ so that $\max\{1 - \alpha, \gamma\}$ is minimal.
It is easy to see that this requires $1 - \alpha = \gamma$, which in turn gives
\begin{equation*}
\alpha = -\frac1{(1 + \beta)\log b}\log\!\left(\frac{1 + \zeta_b(\beta)}{b}\right) .
\end{equation*}
Note that this choice indeed satisfies $\alpha \in {]0, 1[}$, as required in our previous arguments.
Hence, we have to choose $\beta$ in order to minimize
\begin{equation*}
\gamma = 1 + \frac1{(1 + \beta)\log b}\log\!\left(\frac{1 + \zeta_b(\beta)}{b}\right) .
\end{equation*}
Since
\begin{equation*}
\frac{\partial \gamma}{\partial \beta} = \frac{F_b(\beta)}{(1 + \beta)^2 \log b} ,
\end{equation*}
by the previous considerations on $F_b(s)$, we get that $\gamma$ is minimal for $\beta = s_{b,0}$.
Thus, we make this choice for $\beta$, so that $1 - \alpha = \gamma = \eta_{b, 0}$.
Finally, putting together \eqref{equ:N1upperbound} and \eqref{equ:N2upperbound}, we obtain
\begin{equation*}
\#\mathcal{N}_{b,0}(x) < x^{1 - \alpha + o(1)} + x^{\gamma + o(1)} < x^{\eta_{b,0} + o(1)}
\end{equation*}
as $x \to +\infty$.
The proof is complete.

\subsection*{Proof of Theorem~\ref{thm:Nbupper}}

The proof of Theorem~\ref{thm:Nbupper} proceeds similarly to the one of Theorem~\ref{thm:Nb0upper}.
We highlight just the main differences.
First, we shall prove that, for $b \geq 4$, Equation~\eqref{equ:Equb} has a unique positive solution.
For $s \geq 0$, define
\begin{equation*}
G_b(s) := \frac{(1 + s)\zeta_b^\prime(s)}{\zeta_b(s)} - \log\!\left(\frac{\zeta_b(s)}{b}\right).
\end{equation*}
Since $b \geq 4$, we have
\begin{equation}\label{equ:Gb0oo}
G_b(0) = -\log\!\left(\left(1-\frac1{b}\right)(b-1)!^{1/(b-1)}\right) < 0 \quad\text{ and }\quad \lim_{s \,\to\, +\infty} G_b(s) = \log b > 0 .
\end{equation}
Furthermore, a bit of computation shows that
\begin{equation}\label{equ:Gbprime}
G_b^\prime(s) = \frac{(1 + s)(\zeta_b(s)\zeta_b^{\prime\prime}(s) - (\zeta_b^\prime(s))^2)}{(\zeta_b(s))^2} > 0 ,
\end{equation}
for all $s \geq 0$, since \eqref{equ:CS}.
Therefore, by \eqref{equ:Gb0oo} and \eqref{equ:Gbprime}, Equation~\eqref{equ:Equb} has a unique positive solution.
Note also that $G_3(0) > 0$, so that $G_3(s) > 0$ for all $s \geq 0$.

Let $\alpha \in {]0,1[}$ be a constant (depending on $b$) to be determined later, and define $\mathcal{N}_b^{\,\prime}(x)$ as in the proof of Theorem~\ref{thm:Nb0upper}.
Hence, by the previous arguments, the bound \eqref{equ:N1upperbound} holds.

Suppose $n \in \mathcal{N}_{b}^{\,\prime\prime}(x) := \mathcal{N}_{b}(x) \setminus \mathcal{N}_b^{\,\prime}(x)$.
This time, put $N := \lfloor \log n / \log b \rfloor + 1$ (instead of $N := \lfloor \log x / \log b \rfloor + 1$), so that $n$ has exactly $N$ base $b$ digits.
For each $d \in \{1, \dots, b - 1\}$, let $N_d$ be the number of base $b$ digits of $n$ which are equal to $d$.
Note that, since $p_b(n) \mid n$, we have $p_b(n) \neq 0$, that is, all the base $b$ digits of $n$ are nonzero.
Hence, $N_1, \dots, N_{b-1}$ are nonnegative integers such that $N_1 + \cdots + N_{b - 1} = N$.
Furthermore,
\begin{equation*}
p_b(n) = 1^{N_1} \cdots (b-1)^{N_{b-1}} \leq x^{\alpha} < b^{\alpha N} .
\end{equation*}
Let $\beta > 0$ be a constant (depending on $b$) to be determined later.
Summing over all possible values for $N_1, \dots, N_{b - 1}$ and $N$, we get
\begin{align*}
\#\mathcal{N}_b^{\,\prime\prime}(x) &\leq \sum_{N \,=\, 1}^{\lfloor \log x / \log b \rfloor + 1} \sum_{\substack{N_1 + \cdots + N_{b - 1} \,=\, N \\ 1^{N_1} \cdots (b-1)^{N_{b - 1}} \,\leq\, b^{\alpha N}}} \frac{N!}{N_1! \cdots N_{b-1}!} \\
&\leq \sum_{N \,=\, 1}^{\lfloor \log x / \log b \rfloor + 1} \sum_{N_0 + \cdots + N_{b - 1} \,=\, N} \frac{N!}{N_1! \cdots N_{b-1}!} \left(\frac{b^{\alpha N}}{1^{N_1} \cdots (b-1)^{N_{b-1}}}\right)^\beta \\
&= \sum_{N \,=\, 1}^{\lfloor \log x / \log b \rfloor + 1} (b^{\alpha\beta} \zeta_b(\beta))^N \ll (b^{\alpha\beta} \zeta_b(\beta))^{\log x / \log b} ,
\end{align*}
and consequently
\begin{equation}\label{equ:N3upperbound}
\#\mathcal{N}_b^{\,\prime\prime}(x) < x^{\delta + o(1)} ,
\end{equation}
as $x \to +\infty$, where
\begin{equation}\label{equ:deltaalphabeta}
\delta := \alpha \beta + \frac{\log \zeta_b(\beta)}{\log b} .
\end{equation}
At this point, in light of \eqref{equ:N1upperbound} and \eqref{equ:N3upperbound}, we shall choose $\alpha$ and $\beta$ so that $\max\{1 - \alpha, \delta\}$ is minimal.
This requires $1 - \alpha = \delta$, which in turn yields
\begin{equation*}
\alpha = - \frac1{(1 + \beta)\log b}\log\!\left(\frac{\zeta_b(\beta)}{b}\right) .
\end{equation*}
Note that this choice indeed satisfies $\alpha \in {]0, 1[}$, as required in our previous arguments.
Hence, we have to minimize
\begin{equation*}
\delta = 1 + \frac1{(1 + \beta)\log b}\log\!\left(\frac{\zeta_b(\beta)}{b}\right) .
\end{equation*}
We have
\begin{equation*}
\frac{\partial \delta}{\partial \beta} = \frac{G_b(\beta)}{(1 + \beta)^2 \log b} .
\end{equation*}
Hence, by the previous considerations on $G_b(s)$, for $b = 3$ we have to choose $\beta = 0$, while if $b \geq 4$ we have to choose $\beta = s_b$.
Making this choice, we get $1 - \alpha = \delta = \eta_b$.
Finally, putting together \eqref{equ:N1upperbound} and \eqref{equ:N3upperbound}, we obtain
\begin{equation*}
\#\mathcal{N}_b(x) < x^{1 - \alpha + o(1)} + x^{\delta + o(1)} < x^{\eta_b + o(1)}
\end{equation*}
as $x \to +\infty$.
The proof is complete.

\section{Lower bound}

\begin{thm}\label{thm:Nb0lower}
Let $b \geq 3$ be an integer.
We have
\begin{equation}\label{equ:Nb0lower}
\#\mathcal{N}_{b,0}(x) > x^{\rho_{b,0} + o(1)} ,
\end{equation}
as $x \to +\infty$, where
\begin{equation}\label{equ:rhob0}
\rho_{b, 0} := \sup_{\alpha_0, \dots, \alpha_{b - 1}} \frac{\left(\sum_{d \,=\, 1}^{b-1} \alpha_d \right)\log\!\left(\sum_{d \,=\, 1}^{b-1} \alpha_d\right) - \sum_{d \,=\, 1}^{b-1} \alpha_d \log \alpha_d}{\left(1 + \sum_{d \,=\, 1}^{b-1} \alpha_d \right)\log b}
\end{equation}
with $\alpha_0, \dots, \alpha_{b - 1} \geq 0$ satisfying
the conditions
\begin{equation}\label{equ:alphad}
\begin{cases}
\alpha_d = 0 & \text{if $d > 1$ and $p \mid d$, $p \nmid b$ for some prime $p$}, \\
\sum_{d \,=\, 2}^{b-1} \alpha_d \nu_p(d) \leq 1 & \text{for all primes $p \mid b$} ,
\end{cases}
\end{equation}
and with the convention $0 \cdot \log 0 := 0$.
\end{thm}

We remark that if $b$ is a prime number then the bound of Theorem~\ref{thm:Nb0lower} is obvious.
Indeed, the primality of $b$ implies $\alpha_d = 0$ for each $d \in \{2,\dots,b - 1\}$, so that
\begin{equation*}
\rho_{b, 0} = \sup_{\alpha_0, \alpha_1 \geq 0} \frac{(\alpha_0 + \alpha_1)\log(\alpha_0 + \alpha_1) - \alpha_0 \log \alpha_0 - \alpha_1 \log \alpha_1}{(1 + \alpha_0 + \alpha_1)\log b}= \frac{\log 2}{\log b} ,
\end{equation*}
and the bound is
\begin{equation}\label{equ:Nb0lowertrivial}
\#\mathcal{N}_{b, 0}(x) > x^{\log 2 / \log b + o(1)} ,
\end{equation}
as $x \to +\infty$.
However, the bound \eqref{equ:Nb0lowertrivial} follows just by considering that $\mathcal{N}_{b, 0}(x)$ contains all positive integers having their base $b$ digits in $\{0,1\}$.

If $b$ is not a prime number, then Theorem~\ref{thm:Nb0lower} gives a better bound than \eqref{equ:Nb0lowertrivial}.
In particular, for $b = 10$, conditions \eqref{equ:alphad} become
\begin{equation}\label{equ:alphad10}
\begin{cases}
\alpha_3 = \alpha_6 = \alpha_7 = \alpha_9 = 0 , \\
\alpha_2 + 2\alpha_4 + 3\alpha_8 \leq 1 , \\
\alpha_5 \leq 1 ,
\end{cases}
\end{equation}
and the right-hand side of \eqref{equ:rhob0} can be maximized under the constrains given by \eqref{equ:alphad10} using the method of Lagrange multipliers.
This gives $\rho_{10,0} > 0.526$, for the choice 
\begin{equation*}
\alpha_0 = \alpha_1 = 1.331, \quad \alpha_2 = 0.476, \quad \alpha_4 = 0.170, \quad \alpha_5 = 1, \quad \alpha_8 = 0.060\,.
\end{equation*}
Hence, the lower bound for $\#\mathcal{N}_{10,0}(x)$ of Theorem~\ref{thm:base10} follows.

\subsection{Proof of Theorem~\ref{thm:Nb0lower}}
Let us assume $x \geq 1$ sufficiently large, and let $\alpha_0, \dots, \alpha_{b-1} \geq 0$ be constants (depending on $b$) to be determined later.
Define
\begin{equation*}
s := \left\lfloor\frac{\log x}{(1 + \alpha_0 + \cdots + \alpha_{b-1})\log b}\right\rfloor .
\end{equation*}
Also, let $N_d := \lfloor \alpha_d s\rfloor$ for each $d \in \{0, \dots, b - 1\}$, and put $N := N_0 + \cdots + N_{b-1}$.

Now suppose $m$ is a positive integer with at most $N$ base $b$ digits, and such that exactly $N_d$ of its base $b$ digits are equal to $d$, for each $d \in \{1, \dots, b - 1\}$.
Moreover, put $n := b^s m$.
Clearly, $n \leq b^{s + N} \leq x$ and $b^s \mid n$.
Then, imposing the conditions \eqref{equ:alphad}, we get that
\begin{equation*}
p_{b,0}(n) = 1^{N_1} \cdots (b-1)^{N_{b-1}} \mid b^s \mid n ,
\end{equation*}
so that $n \in \mathcal{N}_{b,0}(x)$.
By elementary combinatorics and by using Stirling's formula, the number of possible values for $m$ is
\begin{align*}
\frac{N!}{N_0! \cdots N_{b-1}!} &= \frac{(\lfloor \alpha_0 s\rfloor + \cdots + \lfloor \alpha_{b-1} s\rfloor)!}{\lfloor \alpha_0 s\rfloor! \cdots \lfloor \alpha_{b-1} s\rfloor!} \\
&= \exp\!\left(s\left(\left(\sum_{d \,=\, 1}^{b-1} \alpha_d\right)\log\!\left(\sum_{d \,=\, 1}^{b-1} \alpha_d\right) - \sum_{d \,=\, 1}^{b-1} \alpha_d \log \alpha_d + o(1)\right)\right) ,
\end{align*}
as $s \to +\infty$.
Hence, lower bound \eqref{equ:Nb0lower} follows.
The proof is complete.

\bibliographystyle{amsplain}

\end{document}